\documentclass[10pt]{amsart}
\usepackage{amsmath}
\usepackage{amssymb}

\def\F{{\bf F}}

\def\G{{\bf G}}

\def\ev{{\rm ev}}

\def\PW{{\rm PW}}
\def\pack{{\rm pack}}
\def\cal{\mathcal}

\def\M{{\bf M}}

\def\N{{\mathbb N}}

\def\SG{{\mathfrak S}}

\def\Des{{\rm Des\,}}

\def\Sym{{\bf Sym}}     
\def\FQSym{{\bf FQSym}} 
\def\WQSym{{\bf WQSym}} 
\def\PQSym{{\bf PQSym}} 
\def\CQSym{{\bf CQSym}} 

\def\Nat{{\bf Nat}}

\newtheorem{example}{Example}[section]

\newtheorem{theorem}[example]{Theorem}

\newtheorem{definition}[example]{Definition}
\newtheorem{proposition}[example]{Proposition}

\newtheorem{lemma}[example]{Lemma}

\def\up#1{\raise 1ex\hbox{\footnotesize#1}}

\def\Proof{\noindent \it Proof -- \rm}

\def\qed{\hspace{3.5mm} \hfill \vbox{\hrule height 3pt depth 2 pt width 2mm}
\bigskip}

\def\U0{{\cal U}_0(gl_N)}

\def\End{{\rm End}}
\def\Surj{{\rm Surj}}
\def\<{\langle}
\def\>{\rangle}

\def\QQ{\hbox{\bf Q}}

\def\K{{\mathbb K}}

\def\DessinMatrix#1{\vcenter{\hbox{\makebox[1.7ex]{$#1$}}}}
\def\GenMatrix#1{\vcenter{\halign{&$\DessinMatrix{##}$\cr#1}}\egroup}

\def\setinterlineskip#1{\baselineskip=0pt
  \lineskip=#1 \lineskiplimit=\maxdimen}

\def\matrice{%
  \bgroup
  \let\ =\omit
  \let\\=\cr
  \setinterlineskip{4.0pt}
  \GenMatrix}

\def\DessinsMatrix#1{\vcenter{\hbox{\makebox[1.3ex]{$\scriptstyle#1$}}}}
\def\GensMatrix#1{\vcenter{\halign{&$\DessinsMatrix{##}$\cr#1}}\egroup}

\def\smallmatrice{%
  \bgroup
  \let\ =\omit
  \let\\=\cr
  \setinterlineskip{3.0pt}
  \GensMatrix}

\newlength{\Hackl}
\newcommand{\Hack}{\vrule height \Hackl width 0pt}
\newcommand{\indexmat}%
    {\smallmatrice{\Hack a\\\Hack b\\\Hack c\\\Hack d\\\Hack e\\}}

\newdimen\Squaresize \Squaresize=14pt
\newdimen\Thickness \Thickness=0.5pt

\def\Square#1{\hbox{\vrule width \Thickness
   \vbox to \Squaresize{\hrule height \Thickness\vss
      \hbox to \Squaresize{\hss#1\hss}
   \vss\hrule height\Thickness}
\unskip\vrule width \Thickness}
\kern-\Thickness}

\def\Vsquare#1{\vbox{\Square{$#1$}}\kern-\Thickness}

\newcommand{\bm}[1]{\mbox{\boldmath{$#1$}}}

\title[Endomorphisms of quasi-shuffle algebras]%
{Natural endomorphisms of quasi-shuffle Hopf algebras}

\author[J.-C. Novelli, F. Patras, J.-Y. Thibon]%
{Jean-Christophe Novelli, Fr\'ed\'eric Patras,\\  Jean-Yves Thibon}

\address[Novelli, Thibon]{Institut Gaspard-Monge, Universit\'e Paris-Est, 
5, boulevard Descartes \\Champs-sur-Marne \\77454 Marne-la-Vall\'ee cedex 2 \\
France}
\address[Patras]{Laboratoire de Math\'ematiques J.A. Dieudonn\'e,\\
Universit\'e de Nice - Sophia Antipolis,\\
Parc Valrose, 06108 Nice Cedex 02,\\ France}
\email[J.-C. Novelli]{novelli@univ-mlv.fr}
\email[F. Patras]{patras@unice.fr}
\email[J.-Y. Thibon]{jyt@univ-mlv.fr}

\begin{document}

\begin{abstract}
The Hopf algebra of  word-quasi-symmetric functions ($\WQSym$), a
noncommutative generalization of the Hopf algebra of quasi-symmetric
functions, can be endowed with an internal product that has several
compatibility properties with the other operations on $\WQSym$. This extends
constructions familiar and central in the theory of free Lie algebras,
noncommutative symmetric functions and their various applications fields, and
allows to interpret $\WQSym$ as a convolution algebra of linear endomorphisms
of quasi-shuffle algebras.
We then use this interpretation to study the fine structure of quasi-shuffle
algebras (MZVs, free Rota-Baxter algebras...). In particular, we compute their
Adams operations and prove the existence of generalized Eulerian idempotents,
that is, of a canonical left-inverse to the natural surjection map to their
indecomposables, allowing for the combinatorial construction of free
polynomial generators for these algebras.
\end{abstract}

\maketitle

\section{Introduction}

Quasi-shuffles appeared as early as 1972 in the seminal approach by P. Cartier
to Baxter algebras (now most often called Rota-Baxter algebras)
\cite{cartier}. Their study was revived and intensified during the last 10
years, for a variety of reasons. The first one was the study of MZVs (multiple
zeta values), for example in the works of Hoffman, Minh, Racinet or Zagier,
since quasi-shuffles encode one presentation of their products. Another line
of study, largely motivated by the recent works of Connes and Kreimer on the
structures of quantum field theories, was the revival of the theory of
Rota-Baxter algebras initiated by M. Aguiar, K. Ebrahimi-Fard, L. Guo, and
others. We refer to \cite{hofQs,racinet,cm08,ikz,am,guo,eg08,egk}, also for
further bibliographical and historical references on these subjects.

The present work arose from the project to understand the combinatorial
structure of ``natural'' operations acting on the algebra of MZVs and, more
generally on quasi-shuffle algebras.
It soon became clear to us that the Hopf algebra of word quasi-symmetric
functions ($\WQSym$), was the right setting to perform this analysis and that
many properties of the classical Lie calculus (incorporated in the theory of
free Lie algebras and connected topics) could be translated into this
framework.

This article is a first step in that overall direction. It shows that word
quasi-symmetric functions act naturally on quasi-shuffle algebras and that
some key ingredients of the classical Lie calculus such as Solomon's Eulerian
idempotents can be lifted to remarquable elements in $\WQSym$. In the process,
we show that $\WQSym$ is the proper analogue of the Hopf algebra $\FQSym$ of
free quasi-symmetric functions (also known as the Malvenuto-Reutenauer Hopf
algebra) in the setting of quasi-shuffle algebras. Namely, we prove a
Schur-Weyl duality theorem for quasi-shuffle algebras extending naturally the
classical one (which states that the linear span of permutations is the
commutant of endomorphims of the tensor algebra over a vector space $V$
induced by linear endomorphims of $V$).

The main ingredient of this theory is that the natural extension of the
internal product on the symmetric group algebras to a product on the linear
span of surjections between finite sets, which induces a new product on
$\WQSym$, is a lift in $\WQSym$ of the composition of linear endomorphisms of
quasi-shuffle algebras. This simple observation yields ultimately the correct
answer to the problem of studying the formal algebraic structure of
quasi-shuffle algebras from the Lie calculus point of view. 

\section{Word quasi-symmetric functions}
\label{wqsym}

In this section, we briefly survey the recent theory of noncommutative
quasi-symmetric functions and introduce its fundamental properties and
structures. The reader is referred to \cite{Hiv-adv,NCSF7,HNT} for details and
further information. Let us mention that the theory of word quasi-symmetric
functions is very closely related to the ones of Solomon-Tits algebras and
twisted descents, the development of which was motivated by the geometry of
Coxeter groups, the study of Markov chains on hyperplane arrangements and
Joyal's theory of tensor species. We will not consider these application
fields here and refer to \cite{tits,brown,schocker,ps2006}.

Let us first recall that the Hopf algebra of noncommutative symmetric
functions~\cite{GKLLRT} over an arbitrary field $\K$ of characteristic zero,
denoted here by $\Sym$, is defined as the free associative algebra over an
infinite sequence $(S_n)_{n\ge 1}$, graded by $\deg S_n=n$, and endowed with
the coproduct
\begin{equation}
\Delta S_n=\sum_{k=0}^nS_k\otimes S_{n-k}\quad \text{(where $S_0$=1)}\,.
\end{equation}
It is naturally endowed with an internal product $*$ such
that each homogeneous component $\Sym_n$ gets identified with the (opposite)
Solomon descent algebra of $\SG_n$, the symmetric group of order $n$.
Some bigger Hopf algebras containing $\Sym$ in a natural way are also
endowed with internal products, whose restriction to $\Sym$ coincides
with $*$. An almost tautological example is $\FQSym$, which, being based on
permutations, with the group law of $\SG_n$ as internal product, induces
naturally the product of the descent algebra \cite{NCSF6,MR}.

A less trivial example~\cite{NTtridend} is $\WQSym^*$, the graded dual of
$\WQSym$ (Word Quasi-symmetric functions, the invariants of the
quasi-sym\-metriz\-ing action on words \cite{NCSF7}). It can be shown that
each homogenous component $\WQSym_n^*$ can be endowed with an internal product
(with a very simple combinatorial definition), for which it is anti-isomorphic
with the Solomon-Tits algebra, so that it contains $\Sym_n$ as a
$*$-subalgebra in a non-trivial way. The internal product of $\WQSym^*$ is
itself induced by the one of $\PQSym$ (parking functions), whose restriction
to the Catalan subalgebra $\CQSym$ again contains $\Sym$ in a nontrivial
way~\cite{NT1}.

Let us recall the relevant definitions. We denote by $A=\{a_1<a_2<\dots\}$ an
infinite linearly ordered alphabet and $A^*$ the corresponding set of words.
The \emph{packed word} $u=\pack(w)$ associated with a word $w\in A^*$ is
obtained by the following process. If $b_1<b_2<\dots <b_r$ are the letters
occuring in $w$, $u$ is the image of $w$ by the homomorphism $b_i\mapsto a_i$.
A word $u$ is said to be \emph{packed} if $\pack(u)=u$. We denote by $\PW$ the
set of packed words.
With such a word, we associate the noncommutative polynomial
\begin{equation}
\M_u(A) :=\sum_{\pack(w)=u}w\,.
\end{equation}
For example, restricting $A$ to the first five integers,
\begin{equation}
\begin{split}
\M_{13132}(A)=&\ \ \ 13132 + 14142 + 14143 + 24243 \\
&+ 15152 + 15153 + 25253 + 15154 + 25254 + 35354.
\end{split}
\end{equation}
As for classical symmetric functions, the nature of the ordered alphabet $A$
chosen to define word quasi-symmetric functions $\M_u(A)$ is largely
irrelevant provided it has enough elements. We will therefore often omit the
$A$-dependency and write simply $\M_u$ for $\M_u(A)$, except when we want to
emphasize this dependency (and similarly for the other types of generalized
symmetric functions we will have to deal with).

Under the abelianization
$\chi:\ \K\langle A\rangle\rightarrow\K[A]$, the $\M_u$ are mapped to the
monomial quasi-symmetric functions $M_I$, where $I=(|u|_a)_{a\in A}$ is the
composition (that is, the sequence of integers) associated with the so-called
evaluation vector $\ev(u)$ of $u$ ($\ev(u)_i:=|u|_{a_i}:=|\{j,u_j=a_i\}|$).
Recall, for the sake of completeness, that the $M_I$ are defined, for
$I=(i_1,...,i_k)$, by:
\begin{equation}
M_I:=\sum\limits_{j_1<...<j_k}a_{j_1}^{i_1}...a_{j_k}^{i_k}.
\end{equation}

The polynomials $\M_u$ span a subalgebra of $\K\langle A\rangle$, called
$\WQSym$ for Word Quasi-Symmetric functions~\cite{FH}. This algebra can be
understood alternatively as the algebra of invariants for the noncommutative
version \cite{NCSF7} of Hivert's quasi-symmetrizing action, which is defined
in such a way that two words are in the same $\SG(A)$-orbit (where $\SG(A)$ is
the group of set automorphisms of $A$) iff they have the same packed word. We
refer to \cite{Hiv-adv} for details on the quasi-symmetrizing action.

As for $\Sym$, $\WQSym$ carries naturally a Hopf algebra structure. Its
simplest definition is through the use of two ordered countable alphabets, say
$A=\{a_1<...<a_n<...\}$ and $B:=\{b_1<...<b_n<...\}$.
Let us write $A+B$ for the ordinal sum of $A$ and $B$ (so that for arbitrary
$i,j$, we have $a_i<b_j$). The unique associative algebra map $\mu$ from
$\K\langle A+B\rangle$ to $\K\langle A\rangle \otimes \K\langle B\rangle$
acting as the identity map on $A$ and $B$ induces a map from $\WQSym (A+B)$ to
$\WQSym(A)\otimes\WQSym(B)\cong \WQSym(A)\otimes\WQSym(A)$,
\begin{equation}
\Delta(\M_u(A+B)):=\sum_{\pack(w)=u}w_{|A}\otimes w_{|B},
\end{equation}
which can be shown to define a Hopf algebra structure on $\WQSym$.  Here, for
an arbitrary subset $S$ of $A+B$, $u_{|S}$ stands for the word obtained from
$u$ by erasing all the letters that do not belong to $S$.

The explicit formula for the coproduct $\Delta$ generalizes the usual one for
free quasi-symmetric functions and the Malvenuto-Reutenauer algebra and reads,
for a packed word $u$ on the interval $[1,n]$:
\begin{equation}
\Delta(\M_u) :=
  \sum_{i=0}^n \M_{u_{|[1,i]}}\otimes\M_{\pack(u_{|[i+1,n]}))}.
\end{equation}

Noncommutative symmetric functions (the elements of $\Sym$), although they can
be defined abstractly in terms of a family of algebraically free generators
$S_n$ (see the begining of this section), also does admit a standard
realization in terms of an ordered alphabet $A$ by
\begin{equation}
S_n(A) =\sum_{w\in A^n,\ \Des(w)=\emptyset}w,
\end{equation}
where $\Des(w)=\{i|w_i>w_{i+1}\}$ denotes the descent set of $w$.
Thus, there is  a natural embedding of $\Sym$ into $\WQSym$:
\begin{equation}\label{stdemb}
S_n =\sum_{\Des(u)=\emptyset}\M_u,
\end{equation}
where the summation is implicitely restricted to packed words $u$ of the
suitable length, that is, here, of length $n$ (and similarly in the
forthcoming formulas).
This embedding extends multiplicatively: 
\begin{equation}
S_{n_1}...S_{n_k}=\sum_{\Des(u)\subset\{n_1,...,n_1+...+n_{k-1}\}}\M_u.
\end{equation}
In terms of the realization of noncommutative symmetric functions and word
quasi-symmetric functions over ordered alphabets $A$, these equalities are
indeed equalities: both sides are formal sum of words with certain shapes.

For later use, let us also mention that the last formula implies (by a
standard M\"obius inversion argument that we omit, see \cite{GKLLRT} for
details on the ribbon basis) that the elements
\begin{equation}
\sum_{\Des(u)=\{n_1,...,n_1+...+n_{k-1}\}} \M_u
\end{equation}
belong to $\Sym$. They form a basis, the ribbon basis $R_I$ of $\Sym$:
for $I=(i_1,...,i_k)$,
\begin{equation}
R_I:=\sum_{\Des(u)=\{i_1,...,i_1+...+i_{k-1}\}} \M_u.
\end{equation}

At last, the embedding of $\Sym$ in $\WQSym$ we just considered is a Hopf
embedding: one can either check directly that the image $\bf S$ of the series
$\sum\limits_{n=0}^\infty S_n$ in $\WQSym$ is grouplike
($\Delta({\bf S})={\bf S}\otimes {\bf S}$), or think in terms of ordered
alphabets and notice that the the coproduct of $\WQSym$ given by the ordinal
sum $A+B$ restricts to a coproduct on $\Sym$ that agrees with the one
introduced at the begining of the section.

\section{Extra structures on $\WQSym$}

A natural question arises from our previous account of the theory: does there
exists an internal product on $\WQSym$ extending the one of $\Sym$? The
question will appear later to be closely connected to the problem of using
$\WQSym$ in order to investigate Hopf algebraic properties, very much as
$\Sym$ (and the dual notion of descent algebras) is classically used to
investigate the properties of tensor spaces and connected graded commutative
or cocommutative Hopf algebras.

It turns out that if we want, for example, an interpretation of $\WQSym$
analogous to that of $\FQSym$ as a convolution algebra of endomorphisms of
tensor spaces \cite{Re,MR,NCSF6}, we have to relax the requirement that the
internal product extend the one of $\Sym$. The two internal products will
coincide only on a certain remarquable subalgebra of infinite series (see
Section~\ref{carac}).
Moreover, the construction does not work with the standard
embedding~(\ref{stdemb}), and therefore we will also have to relax the
condition that the embedding of $\Sym$ in $\WQSym$ is compatible with
realizations in terms of ordered alphabets. This results into a new picture of
the relations between $\Sym$ and $\WQSym$, where one has to map the complete
symmetric functions as follows
\begin{equation}
\label{goodemb}
S_n \mapsto \M_{12\ldots n}\,.
\end{equation}
Requiring the map to be multiplicative defines a new morphism of algebras -for
an element $T\in\Sym$, we will write $\hat T$ its image in $\WQSym$.

\begin{proposition} 
This map is still a Hopf algebra embedding. Its action on the monomial basis
generated by the $S_n$ and on the ribbon basis is given respectively by: For
$I=(i_1,...,i_k)$ with $i_1+...+i_k=n$,
\begin{equation}
\label{rubs0}
\hat S^I := \hat S_{i_1}...\hat S_{i_k}
 = \sum_{\Des(u)\supseteq[n]-\{i_k,i_k+i_{k-1},...,i_k+...+i_1\}}
       \M_{\check u}\,,
\end{equation}

\begin{equation}
\label{rubs}
\hat R_I
 = \sum_{\Des(u)=[n]-\{i_k,i_k+i_{k-1},...,i_k+...+i_1\}}\M_{\check u}\,,
\end{equation}
where $w=a_{j_1}\cdots a_{j_n}\mapsto\check{w}=a_{j_n}\cdots a_{j_1}$ is the
anti-automorphism of the free associative algebra reversing the words.
\end{proposition}

The first assertion follows from the observation that two series $\sum_nS_n$
and $\sum_n\M_{12\ldots n}$ are grouplike and generate two free associative
algebras (respectively $\Sym$ and a subalgebra of $\WQSym$ isomorphic to
$\Sym$).

Let us compute the action on the monomial basis (the action on the ribbon
basis follows by M\"obius inversion since
$\hat S^I=\sum\limits_{J\subset I}\hat R_J$).
From $\hat S_n=\sum\limits_{i_1<...<i_n}a_{i_1}...a_{i_n}$, we get:
\begin{equation}
\hat S_n=\sum_{\Des(u)=[n-1]}\M_{\check u}.
\end{equation}
The same principle applies in general (if a word of length $n$ is strictly
increasing in position $i$, then the reverse word has a descent in position
$n-i$) and we get:
\begin{equation}
\begin{split}
\hat S^I
&=\sum\limits_{j_1^1<\dots<j_{i_1}^1,\dots,j_1^k<\dots<j_{i_k}^k}
     a_{j_1^1}...a_{j_{i_k}^k}
\\
&=\sum\limits_{\Des(u)\supseteq[n]-\{i_k,i_k+i_{k-1},\dots,i_k+...+i_1\}}
   \M_{\check u},
\end{split}
\end{equation}
\begin{equation}
\end{equation}
from which the Proposition follows.

Moving beyond the relation to $\Sym$, recall that, in general,
\begin{equation}
\M_u\M_v =\sum_{w=u'v',\pack(u')=u,\ \pack(v')=v}\M_w\, .
\end{equation}
An interesting feature of  $\WQSym$ is the quasi-shuffle nature of this
product law. This feature explains many of its universal properties with
respect to quasi-shuffle algebras.

To understand it, first notice that packed words $u$ over the integers (recall
that the alphabet can be chosen arbitrarily provided it is ``big enough'') can
be interpreted as surjective maps
\begin{equation}
u:\ [n]\longrightarrow [k]\,,\quad (k=\max(u))\, \quad u(i):=u_i
\end{equation} 
or, equivalently, as ordered partitions (set compositions) of $[n]$: 
let us write $\overline u$ for $(u^{-1}(1),...,u^{-1}(k))$.
We will use freely these two interpretations of packed words from now on to
handle computations in $\WQSym$ as computations involving surjective morphisms
or ordered partitions.

Now, let us define recursively the quasi-shuffle product
$S\uplus T=(S_1,...,S_k) \uplus (T_1,...,T_l)$ of two sequences of sets $U$
and $V$ (which is a formal sum of sequences of sets) by:
\begin{equation}
U\uplus V:=(U_1,U'\uplus V)+(V_1,U\uplus V')+(U_1\cup V_1,U'\uplus V'),
\end{equation}
where $U':=(U_2,...,U_k)$ and $V':=(V_2,...,V_l)$.
Then, if $u$ (resp. $v$) encodes the set composition $U$ (resp. $V$), in
$\WQSym$:
\begin{equation}
\M_u\M_v=\M_t,
\end{equation}
where the set of packed words $t$ encodes $T= U\uplus V[n]$, $n$ is the length
of $u$, and for an arbitrary sequence $S$ of subsets of the integers,
$S[p]:=(S_1+p,...,S_k+p)$.
We will write abusively $\uplus$ for the operation on the linear span of
packed words induced by the ``shifted shuffle product''
$\overline u\uplus \overline v[n]$ so that, with our previous conventions,
$t=u\uplus v$ and (with a self-explanatory notation for $\M_{u\uplus v}$)
$\M_u\M_v=\M_{u\uplus v}$.

To conclude this section, let us point out that a candidate for the internal
product we were looking for is  easily described using the interpretation of
packed words as surjections:
\begin{definition}
The internal product of $\WQSym$ is defined in the $\M$-basis by
\begin{equation}
\label{defint}
\M_u*\M_v 
  = \M_{v\circ u}\quad \text{whenever $l(v)=\max(u)$ and $0$ otherwise.}
\end{equation} 
\end{definition}
The following sections show that this product has the expected properties with
respect to the other structures of $\WQSym$ and with respect to arbitrary
quasi-shuffle algebras.

\section{A relation between internal and external products}

In $\Sym$, there is a fundamental compatibility relation between the internal
product, the usual product and the coproduct.
It is called the splitting formula \cite{GKLLRT}, and is essentially a
Hopf-algebraic interpretation of the noncommutative Mackey formula discovered
by Solomon \cite{So}. It can be extended to $\FQSym$, with certain
restrictions \cite{NCSF7}. The key ingredient for doing this is an expression
of the product of $\FQSym$ in terms of shifted concatenation and internal
product with an element of $\Sym$. This can again be done here.

The natural notion of shifted concatenation in $\WQSym$ is not the same as in
$\FQSym$: indeed, if $u$ and $v$ are packed words, one would like that
$u\bullet v$ be a packed word.  The correct way to do this is to shift the
letters of $v$ by the maximum of $u$:
\begin{equation}
u\bullet v = u\cdot v[k]\,,\quad\text{where $k=\max(u)$}\,.
\end{equation}
For example, $11\bullet 21=1132$. We consistently set
\begin{equation}
\M_u\bullet \M_v =\M_{u\bullet v}.
\end{equation}

\begin{lemma}
We have the distributivity property: 
\begin{equation}
 (\M_u\bullet\M_t)\ast (\M_v\bullet \M_w)
=(\M_u\ast \M_v)\bullet (\M_t\ast \M_w),
\end{equation}
whenever $l(v)=\max(u), l(w)=\max(t)$.
\end{lemma}

We also have the following crucial lemma (compare with~\cite[Eq. (2)]{NCSF7}):
\begin{lemma}
\label{Lcrucial}
Let $u_1,\ldots,u_r$ be packed words, and define a composition
$I=(i_1,\ldots,i_r)$ by $i_k=\max(u_k)$. Then, if $\Sym$ is embedded in
$\WQSym$ by means of~(\ref{goodemb}), and the internal product is defined
by~(\ref{defint}),
\begin{equation}
\label{crucial}
\M_{u_1}\M_{u_2}\cdots \M_{u_r} =
(\M_{u_1}\bullet\M_{u_2}\bullet \cdots\bullet  M_{u_r} ) * S^I\,.
\end{equation}
\end{lemma}

For example,
\begin{equation}
\M_{11}\M_{21}= \M_{1132}+ \M_{1121}+\M_{2231}+\M_{2221}+\M_{3321}
\end{equation}
is obtained from
\begin{equation}
\M_{11}\bullet \M_{21}=\M_{1132}
\end{equation}
by internal product on the right by
\begin{equation}
S^{12}=S_1S_2=\M_1\M_{12}=\M_{123}+\M_{112}+\M_{213}+\M_{212}+\M_{312}\,.
\end{equation}

\Proof
The Lemma is most easily proven by switching to the langage of surjections.
Let us notice first that, by construction of the shifted quasi-shuffle product
of set compositions, $u_1\uplus  u_2$ is the formal sum of all surjections
from $[l(u_1)+l(u_2)]$ to $[i_1+i_2-p]$ , where $p$ runs from $0$ to $\inf
(i_1,i_2)$, that can be obtained by composition of $u\bullet v$ with a
surjective map $\phi$ from $[i_1+i_2]$ to $[i_1+i_2-p]$ such that $\phi$ is
(strictly) increasing on $[i_1]$ and $\{i_1+1,...,i_1+i_2\}$.

Let us write $\gamma_{i_1,i_2}$ for the formal sum of these surjections with
domain $[i_1+i_2]$ and codomain $[i_1+i_2-p]$, $p=0,\dots,\inf(i_1,i_2)$. We
get, as a particular case: $1_{i_1}\uplus 1_{i_2}=\gamma_{i_1,i_2}$, where
$1_n$ stands for the identity in $\SG_n$, the symmetric group of rank $n$.
In general, we have therefore:
\begin{equation}
 \M_{u_1}\M_{u_2}=(\M_{u_1} \bullet \M_{u_2})\ast
   (\M_{1\dots i_1}\M_{1\dots i_2})
=(\M_{u_1} \bullet \M_{u_2})\ast S^{i_1,i_2},
\end{equation}
with the notation
\begin{equation}
S^{i_1\dots i_k}
 := S_{i_1}\dots S_{i_k}
  = \M_{1\dots i_1} \dots \M_{1\dots i_k}.
\end{equation}
The same reasoning applies to an arbitrary number of factors.
\qed

\section{Operations on quasi-shuffle algebras}

Let us recall first the definition of the quasi-shuffle algebra $QS(A)$ on a
commutative algebra $A$ (without a unit and over $\K$).
The underlying vector space is the tensor algebra over $A$:
$QS(A)=\bigoplus_nA^{\otimes n}$, where $A^{\otimes 0}:=\K$. The product is
defined recursively by:
\begin{equation}
\begin{split}
(a_1\otimes ...\otimes a_n)\uplus (b_1\otimes ...\otimes b_m)
=&a_1\otimes ((a_2\otimes ...\otimes a_n)\uplus (b_1\otimes ...\otimes b_m))\\
&+b_1\otimes ((a_1\otimes ...\otimes a_n)\uplus (b_2\otimes ...\otimes b_m))\\
&+a_1b_1\otimes ((a_2\otimes ...\otimes a_n)\uplus (b_2\otimes ...\otimes b_m)).
\end{split}
\end{equation}

\begin{lemma}
The quasi-shuffle algebra is a right module over $\WQSym$ equipped with the internal product. The action is defined by:
\begin{equation}
(a_1\otimes ...\otimes a_n)\M_u
:=\delta_m^n b_1\otimes ...\otimes b_k,\ \ b_i:=\prod_{u(j)=i}a_j,
\end{equation}
where we used the surjection interpretation of packed words, $u$ is a
surjective map from $[m]$ to $[k]$, and $\delta_m^n$ is the Kronecker symbol
($\delta_m^n=1$ if $m=n$ and $=0$ else).
\end{lemma}

This right-module structure allows to rewrite the definition of the
quasi-shuffle product as (see, {\it e.g.},~\cite{cartier}, also for a proof
that $QS(A)$ is actually a commutative algebra):
\begin{equation}
(a_1\otimes ...\otimes a_n)\uplus (b_1\otimes ...\otimes b_m)
=(a_1\otimes ...\otimes a_n\otimes b_1\otimes ...\otimes b_m)S^{n,m}.
\end{equation}
Interpreting $S^{n,m}$ as an element of $\FQSym$ instead of $\WQSym$ and using
the standard right action of permutations on tensors, we would get the ordinary
shuffle product.

\begin{lemma}
The quasi-shuffle algebra is endowed with a Hopf algebra structure by the
deconcatenation coproduct
\begin{equation}
\Delta(a_1\otimes ...\otimes a_n)
:=\sum_{i=0}^n (a_1\otimes ...\otimes a_i)\otimes
               (a_{i+1}\otimes ...\otimes a_n):
\end{equation}
\begin{equation}
 \Delta((a_1\otimes ...\otimes a_n)\uplus (b_1\otimes ...\otimes b_m))
=\Delta(a_1\otimes ...\otimes a_n)(\uplus\otimes\uplus)
  \Delta(b_1\otimes ...\otimes b_m).
\end{equation}
\end{lemma}

This Lemma, due to Hoffman \cite{hofQs}, amounts to checking that both sides
of this last identity are equal to:
\begin{equation}
\sum_{i\leq n,j\leq m}((a_1\otimes ...\otimes a_i)\uplus
                       (b_1\otimes ...\otimes b_j))
\otimes ((a_{i+1}\otimes ...\otimes a_n)\uplus
         (b_{j+1}\otimes ...\otimes b_m)),
\end{equation}
which follows immediately from the definition of $\Delta$ and $\uplus$.

\begin{proposition}
\label{compatibilite}
The right module structure of $QS(A)$ over $\WQSym$ is compatible with the
outer product (\emph{i.e.}, the usual graded product of $\WQSym$, induced by
concatenation of words), in the sense that this product coincides with the
convolution product $\star$ in $\End(QS(A))$ induced by the Hopf algebra
structure of $QS(A)$.
\end{proposition}

Indeed, by definition of the convolution product of Hopf algebras linear
endomorphisms, we have, for $u$ and $v$ surjections from $[n]$ (resp. $[m]$)
to $[p]$ (resp. $[q]$):
\begin{equation}
\begin{split}
(a_1\otimes ...\otimes a_{n+m})(\M_u\star \M_v)
&=((a_1\otimes ...\otimes a_n)\M_u)\uplus
   (a_{n+1}\otimes ...\otimes a_{n+m})\M_v\\
&=(((a_1\otimes ...\otimes a_n)\M_u)\otimes 
   ((a_{n+1}\otimes ...\otimes a_{n+m})\M_v))S^{p,q}\\
&=(a_1\otimes ...\otimes a_{n+m})(\M_u\bullet \M_v)S^{p,q}\\
&=(a_1\otimes ...\otimes a_{n+m})\M_u\M_v\,.
\end{split}
\end{equation}
by Lemma~\ref{Lcrucial}, or, since the identity does not depend on
$a_1,...,a_{n+m}$:
\begin{equation}
\M_u\star \M_v=\M_u\M_v.
\end{equation}

Let us formalize, for further use, our last observation on the dependency on
$a_1,...,a_{n+m}$ into a general recognition principle that will prove useful
to deduce properties in $\WQSym$ from its action on quasi-shuffle algebras.

\begin{lemma}
\label{recognition}
Let $f$ and $g$ be two elements in $\WQSym$ and let us assume that, for an
arbitrary commutative algebra $A$ and arbitrary $a_1,...,a_n,...$ in $A$,
$(a_1\otimes ...\otimes a_n)f=(a_1\otimes \cdots\otimes a_n)g$ for all $n$.
Then, $f=g$.
\end{lemma}

The Lemma follows, \emph{e.g.},  by letting $a_1,...,a_n,...$ run over an
infinite ordered alphabet and letting $A$ be the free commutative algebra over
this alphabet.

\section{Nonlinear Schur-Weyl duality}

The same kind of argument can actually be used to characterize $\WQSym$ as a
universal endomorphism algebra in the same way as $\FQSym$ is a universal
endomorphism algebra according to the classical Schur-Weyl duality.
Recall the latter: for an arbitrary vector space $V$, let us write
$T(V):=\bigoplus\limits_{n\in\N}T_n(V):=\bigoplus\limits_{n\in\N}V^{\otimes
n}$. Linear morphisms between vector spaces $f:V\longmapsto W$ induce maps
$T(f):T(V)\longmapsto T(W)$ compatible with the graduation ($T_n(f)$ maps
$T_n(V) $ to $T_n(W)$: $T_n(f)(v_1\otimes ...\otimes v_n):=f(v_1)\otimes
...\otimes f(v_n)$). In categorical langage, Schur-Weyl duality  characterizes
natural transformations of the functor $T$ (or, equivalently, of the
subfunctors $T_n$) from vector spaces to graded vector spaces and reads: the
only family of maps $\mu_V: T_n(V)\longmapsto T_n(V)$ (where $V$ runs over
vector spaces over $\K$) such that, for any map $f$ as above,
\begin{equation}
T_n(f)\circ \mu_V=\mu_W\circ T_n(f),
\end{equation}
are linear combination of permutations: $\mu_V\in \K [S_n]$ (the converse
statement is obvious: permutations and linear combinations of them acting on
tensors always satisfy this equation).

We consider here the corresponding nonlinear problem and characterize natural
transformations of the functor
\begin{equation}
T(A):=\bigoplus\limits_{n\in\N}T_n(A):=\bigoplus\limits_{n\in\N}A^{\otimes n}
\end{equation}
viewed now as a functor from commutative algebras without a unit to vector
spaces. Concretely, we look for families of linear maps $\mu_A$ from $T_n(A)$
to $ T_m(A)$, where $A$ runs over commutative algebras without a unit and $m$
and $n$ are arbitrary integers such that, for any map $f$ of algebras from $A$
to $B$,
\begin{equation}
\label{natu}
T_m(f)\circ \mu_A=\mu_B\circ T_n(f).
\end{equation}
Let us say that such a family $\mu_A$ satisfies nonlinear Schur-Weyl duality
(with parameters $n,m$).
The purpose of the section is to prove:

\begin{proposition}
Let $\Nat$ be the vector space spanned by families of linear maps that satisfy
the nonlinear Schur-Weyl duality. Then $\Nat$ is canonically isomorphic to
$\WQSym$.

Equivalently, the vector space $\Nat_{n,m}$ of families of linear maps that
satisfy non linear Schur-Weyl duality with parameters $n,m$ is canonically
isomorphic to the linear span of surjections from $[n]$ to $[m]$.
\end{proposition}

The results in the previous section imply that $\WQSym$ is canonically
embedded in $\Nat$. Let us show now that the converse property holds. We write
$\QQ[x_1,...,x_n]^+$ for the vector space of polynomials in the variables
$x_1,...,x_n$ without constant term and notice that, for an arbitrary family
$a_1,...,a_n$ of elements of a commutative algebra $A$, the map $f(x_i):=a_i$
extends uniquely to an algebra map from $\QQ[x_1,...,x_n]^+$ to $A$. In
particular, if $\mu_A$ is a family of linear maps that satisfy the nonlinear
Schur-Weyl duality, we have:
\begin{equation}
\mu_A(a_1\otimes ...\otimes a_n)=\mu_A\circ T(f)(x_1\otimes ...\otimes x_n)
= T(f)(\mu_{\QQ[x_1,...,x_n]^+}(x_1\otimes ...\otimes x_n)),
\end{equation}
so that the knowledge of $\mu_{\QQ[x_1,...,x_n]^+}(x_1\otimes ...\otimes x_n)$
determines entirely the other maps $\mu_A$.

Let $\mu_A\in \Nat_{n,m}$. Then, $\mu_{\QQ[x_1,...,x_n]^+}(x_1\otimes
...\otimes x_n)\in (\QQ[x_1,...,x_n]^+)^{\otimes m}$. 
The latter vector space has a basis $\mathcal B$ whose elements are the
tensors ${\mathbf p}=p_1\otimes ...\otimes p_m$, where the $p_i$s run over all
the nontrivial monomials in the $x_i$s (for example $x_1^2x_3\otimes
x_2^2\otimes x_1x_5$ is a basis element for $n=5$ and $m=3$): 
$\mu_{\QQ[x_1,...,x_n]^+}(x_1\otimes ...\otimes x_n)$ can therefore be written
uniquely a linear  combination of these basis elements. 

Now, the commutation property (\ref{natu}) implies that
\begin{equation}
Y:=\mu_{\QQ[x_1,...,x_n]^+}(x_1\otimes ...\otimes x_n)
  \in (\QQ[x_1,...,x_n]^+)^{\otimes m}
\end{equation}
must be linear in each variable $x_i$ (take $f$ such that $x_i\mapsto ax_i$
and $x_j\mapsto x_j$ for $j\not = i$).

This implies in particular that $\Nat_{n,m}=0$ when $m>n$. For $m=n$, $Y$ must
be a linear combination of permutations
$x_{\sigma(1)}\otimes\cdots\otimes x_{\sigma(n)}$, and for $m<n$, 
\begin{equation}
Y=\sum_{{\mathbf p}\in {\mathcal B},deg({\mathbf p})=n}
  \lambda_{\mathbf p}{\mathbf p}
\end{equation}
with $p_1\dots p_n= x_1\dots x_n$, which implies that
\begin{equation}
\mu_{\QQ[x_1,\dots,x_n]^+}(x_1\otimes \dots\otimes x_n)
=\sum_{f\in \Surj(n,m)}
\lambda_f(\prod\limits_{f(i)=1}x_i\otimes \dots\otimes
          \prod\limits_{f(i)=m}x_i).
\end{equation}
Thus, $\mu_A$ can necessarily be written as a linear combination of (maps
induced by) surjections from $n$ to $m$.
\qed

\section{The Characteristic subalgebra of $\WQSym$}
\label{carac}

The existence of two algebra maps from $\WQSym$ (equipped with the product of
word quasi-symmetric functions and the internal product) to $\End(QS(A))$
(equipped with the convolution product and the composition product) extends a
classical result. There are indeed two analogous maps from $\FQSym$ to the
endomorphism algebra of the tensor algebra over an alphabet $X$ \cite{MR}:
this corresponds roughly to the case where one considers $QS(A)$ with $A$ the
linear span of $X$ equipped with the null product and can be understood as a
particular case of the constructions we are interested in here.

In the ``classical'' situation, it is however well-known that, from the Hopf
algebraic point of view, most relevant informations are contained in a very
small convolution subalgebra of $\FQSym$, namely the one generated by the
identity of the algebra \cite{Pat92,Pat93}. In the present section, we
investigate the structure of the corresponding subalgebra of $\WQSym$ and
deduce from this study that many essential objects in Lie theory (Solomon's
idempotents...) have a quasi-shuffle analogue in $\WQSym$. Most results in
this section are direct applications of \cite[Chap. 1]{Pat92} (published in
\cite{Pat93}), to which we refer for proofs and details.

We make implicitely use of the recognition principle (Lemma \ref{recognition})
to deduce these results from the existence of an action of $\WQSym$ of
$QS(A)$.
To deal with formal power series in $\WQSym$, we consider the usual topology
(the one associated to the graduation induced by word length, that is, the one
for which words with large lengths are close to 0), and write $\hat\WQSym$ for
the corresponding completion.

\begin{lemma}
The $k$-th characteristic endormorphism (or Adams operation) of $QS(A)$ is the
following element of $\hat\WQSym$, identified with the $k$-th convolution
power of the identity map:
\begin{equation}
\Psi^k:=I^{k},\ I:=\hat\sigma_1:=\sum_{n\geq 0}\M_{1,...,n}.
\end{equation}
The characteristic endomorphisms satisfy:
\begin{itemize}
\item $\Psi^k$ is an algebra endomorphism of $QS(A)$,
\item $\Psi^k \Psi^l=\Psi^{k+l}$,
\item $\Psi^k\ast \Psi^l=\Psi^{kl}$.
\end{itemize}
\end{lemma}

See \cite[Prop. 1.4, Prop. 1.3]{Pat93}. To deduce the identity  $\Psi^k
\Psi^l=\Psi^{k+l}$, we use the fact that the product in $\WQSym$ maps to the
convolution product in an arbitrary $QS(A)$.

Many important structure results that hold for graded Hopf algebras
\cite{Pat94} do not hold for quasi-shuffle bialgebras -that are not graded but
only filtered: the product in $QS(A)$ maps $A^{\otimes n}\otimes A^{\otimes
m}$ to $\bigoplus\limits_{p\leq n+m}A^{\otimes p}$, whereas the coproduct
respects the graduation and maps $A^{\otimes p}$ to
$\bigoplus\limits_{n+m=p}A^{\otimes n}\otimes A^{\otimes m}$. However, some
properties of graded Hopf algebras hold for the quasi-shuffle algebras:
\begin{itemize}
\item Let $f,g$ two linear endomorphisms of $QS(A)$ that vanish on
$\bigoplus\limits_{p\leq n}A^{\otimes p}$, resp. $\bigoplus\limits_{p\leq
m}A^{\otimes p}$, then, since the coproduct preserves the graduation,
$f\star g$ vanishes on $\bigoplus\limits_{p\leq n+m+1}A^{\otimes p}$ (we say
that $f,g,f\star g$ are respectively $n,m,n+m+1$-connected).
\item  The element $I$ is invertible in $\hat\WQSym$ -this follows from
\begin{equation}
I^{-1}=(M_0+\sum\limits_{n>0}M_{1...n})^{-1}
=\sum\limits_{k\geq 0}(-1)^k(\sum\limits_{n>0}M_{1...n})^k,
\end{equation}
since $\sum\limits_{n>0}M_{1...n}$ is $0$-connected.
\end{itemize}

\begin{definition}
We call Characteristic subalgebra of $\hat\WQSym$ and write $\mathbf{Car}$ for
the convolution subalgebra of $\hat\WQSym$ generated by $I$.
\end{definition}

\begin{lemma}
The representation of $\mathbf{Car}$ on
$QS^{(n)}(A):=\bigoplus\limits_{p\leq n}A^{\otimes p}$ is unipotent of rank
$n+1$. That is, for any 0-connected element $f$ in $\mathbf{Car}$, $f^{n+1}$
acts on $QS^{(n)}(A)$ as the null operation.
\end{lemma}

This follows from the previously established properties.

\begin{proposition}
The action of $\Psi^k$ on $QS^{(n)}(A)$ is polynomial in $k$:
\begin{equation}
\Psi^k=\sum\limits_{i=0}^nk^ie_i,
\end{equation}
where the ``quasi-Eulerian idempotents'' are given by
\begin{equation}
e_i={\displaystyle \frac{\log(I)^i}{i!}}\,.
\end{equation}
\end{proposition}

This follows from the convolution identity
$\Psi^k=I^k=\exp(\log (I^k))=\exp(k\log(I))$ since
$\log(I)=\sum\limits_{n\geq 1}(-1)^{n+1}\frac{(I-M_0)^n}{n}$ is $0$-connected.
One can make explicit the formula for the $e_i$s using the Stirling
coefficients of the first kind (Fla I.4.5 in \cite{Pat92}).
The $e_i$ are orthogonal idempotents (this follows as in the usual case from
$\Psi^k \Psi^l=\Psi^{k+l}$: the proof of \cite[Prop.I,4,8]{Pat92}
\cite[Prop.3.4]{Pat93} applies). When $A$ is the linear span of an alphabet
equipped with the null product, we recover Solomon's Eulerian idempotents.

Using the fact that
\begin{equation}
\Psi^k =\hat\sigma_1^k
\end{equation}
and computing in $\Sym$, we obtain from (\ref{rubs}) the following expression.
\begin{proposition}
We have:
\begin{equation}
e_1 = \sum_{n\ge 1}
      \frac1n\sum_{I\vDash n}
      \frac{(-1)^{l(I)-1}}{{n-1\choose l(I)-1}}
          \sum_{\Des(u)=[n]-\{i_{l(I)},\dots,i_{l(I)}+\dots+i_1\}}
                \M_{\check u}\,,
\end{equation}
where $I\vDash n$ means that $I=(i_1,...,i_{l(I)})$ is a composition of $n$
($i_1+...+i_{l(I)}=n$).
\end{proposition}

For example, up to degree 3,
\begin{equation}
\begin{split}
e_1 = &\M_1\\
&+\frac12(\M_{12}-\M_{11} -\M_{21})\\
&+\frac16(2\M_{123}- \M_{122}- \M_{112}+2 \M_{111}- \M_{231}- \M_{132}\\
&+2\M_{221}-\M_{121}-\M_{213}-\M_{212}+2 \M_{211}+2 \M_{321}- \M_{132})
 +\cdots
\end{split}
\end{equation}

\begin{equation}
\begin{split}
e_2 =& \frac12 (\M_{12}+  \M_{11} +\M_{21} )\\
&+\frac12( \M_{123}-  \M_{111} - \M_{221} - \M_{211} - \M_{321})+\cdots
\end{split}
\end{equation}

\begin{equation}
e_3 = \frac16\sum_{|u|=3}\M_{u}+\cdots
\end{equation}

\section{The case of quasi-symmetric functions}

The fundamental example of a quasi-shuffle algebra is $QSym$, the Hopf algebra
of quasi-symmetric functions (it is the quasi-shuffle algebra over the algebra
of polynomials in one variable, or in additive notation, over the nonnegative
integers \cite{TU}).
Thus, it is of some interest to have a closer look at the right action of
$\WQSym$ on $QSym$ as defined in the foregoing section. Because of the duality
between $QSym$ and $\Sym$, this will also result into a refined understanding
of the links between word quasi-symmetric functions and noncommutative
symmetric functions from a Lie theoretic point of view.

The basis which realizes $QSym$ as a quasi-shuffle algebra is the
quasi-monomial basis $M_I$ whose definition was recalled in
Section~\ref{wqsym}:
\begin{equation}
M_I M_J = M_{I\uplus J} :=\sum_{K}(K|I\uplus J) M_K\,,
\end{equation} 
where, for two compositions $K$ and $L$, $(K|L):=\delta_K^L$.
Let us denote by a $*$ the right action of $\WQSym$: for $I=(i_1,...,i_k)$ and $u$ a packed word of length $k$,
\begin{equation}
M_I*\M_u = M_J\,, \quad j_r=\sum_{u(s)=r}i_s.
\end{equation} 
For example,
\begin{equation}
M_{21322}*\M_{12121}=M_{2+3+2,1+2}=M_{73}\,.
\end{equation}
Then, the fundamental compatibility formula (Prop.~\ref{compatibilite})
can be rewitten as
\begin{equation}
M_I * (\M_u\M_v) = \mu[\Delta M_I *_2 (\M_u\otimes \M_v)]
\end{equation}
where $\mu$ is the multiplication map.
This mirrors the splitting formula for the internal product of $\FQSym$
\begin{equation}\label{splitF}
(\F_\sigma\F_\tau) * S^I = \mu[(\F_\sigma\otimes \F_\tau) *_2 \Delta S^I]
\end{equation}
which can be extended to any number of factors on the left, and $S^I$ be
replaced by an arbitrary noncommutative symmetric function. Similarly, we can
write, for any $F\in QSym$ and $\G_1,\ldots, \G_r\in\WQSym$, 
\begin{equation}
F * (\G_1\G_2\cdots \G_r)
=\mu_r[\Delta^r F *_r (\G_1\otimes\G_2\otimes\cdots\otimes \G_r)],
\end{equation}
where $\mu_r$, $*_r$ and $\Delta^r$ stand for $r$-fold iterations of the
corresponding product and coproduct maps.
One may ask whether the formula remains valid for bigger quotients of $\WQSym$
(recall that $QSym$ is its commutative image). It appears that one must have
$\M_u\equiv\M_v$ for $\ev(u)=\ev(v)$ except when $u$ and $v$ are formed of an
equal number of 1s and 2s, in which case other choices are allowed.
Thus, $QSym$ is essentially the only interesting quotient.

The commutative image map can be expressed by means of the action $*$. Recall
that $M_{1^n}=e_n$ is the $n$-th elementary symmetric function. For a packed
word $u$ of length $n$,
\begin{equation}
M_{1^n}*\M_u = M_{\ev(u)}
\end{equation}
so that for any $\G\in\WQSym$, its commutative image $G$ is
\begin{equation}
G=\lambda_1*\G \quad \left(\lambda_1:=\sum_{n\ge 0}M_{1^n}\right).
\end{equation}

The Adams operations $\Psi^k$ of $QSym$ defined above are
\begin{equation}
\Psi^k(F) = \mu_k\circ \Delta^k (F) =: F(kX)
\end{equation}
where the $\lambda$-ring theoretical notation $F(kX)$ is motivated by the
observation that the coproduct of $QSym$ can be defined by means of ordinal
sums of alphabets ($\Delta F(X)=F(X+Y)$). By the previous theory, we have
\begin{equation}
F(kX)= F * \hat\sigma_1^k = F*\left(\sum_{n\ge 0}\M_{12...n}\right)^k\,.
\end{equation}
For example, with $k=2$, we can easily compute the first terms by hand
\begin{equation}
\begin{split}
\hat\sigma_1^2= &1 + 2 \M_1 + \M_1\M_1 + 2\M_{12}\\
&+ \M_1\M_{12}+\M_{12}\M_1 + 2\M_{123} + ...\\
=& 1 + 2\M_1\\
&+ 3\M_{12}+\M_{21}+\M_{11}\\
&+ 4\M_{123}+ M_{112}+\M_{213}+\M_{312} +\M_{121}\\
& \ + \M_{122}+ \M_{132}+ \M_{231}+...
\end{split}
\end{equation}
so that $\Psi^2(M_n)=2M_n$, and
$\Psi^2(M_{ij})=3M_{ij}+M_{ji}+M_{i+j}$, {\it etc.},
which agrees indeed with the direct computation $\Psi^k(M_I)=M_I(kX)$,
as $M_{ij}(X+Y)=M_{ij}(X)+M_i(X)M_j(Y)+M_{ij}(Y)$ which
for $Y=X$ does yield $3M_{ij}+M_{ji}+M_{i+j}$.

Note that we can have sligthly more general operators by introducing
extra parameters, \emph{e.g.}, 
\begin{equation}
M_I * \hat\sigma_t = t^{\ell(I)}M_I.
\end{equation}
Then, we would have deformations of the Adams operations, like
\begin{equation}
M_{ij} * (\hat\sigma_x\hat\sigma_y)=(x^2+xy+y^2)M_{ij}+xy(M_{ji}+M_{i+j})\,.
\end{equation}

Let us take now advantage of the duality between $\Sym$ and $QSym$
\cite{GKLLRT}. Take care that the following duality results are specific to
$QSym$ and $\Sym$ and would not hold for arbitrary quasi-shuffle algebras -in
particular, there is no such direct link in general between the quasi-Eulerian
idempotents acting on $QS(A)$ and the usual Eulerian idempotents as the one
described below.

We write $\psi^k$ for the adjoint of $\Psi^k$, acting on $\Sym$. Since the
product and coproduct on $\Sym$ are dual to the ones on $QSym$, we have again,
on $\Sym$, $\psi^k := \mu_k\circ\Delta^k$, where now $\mu_k$ is the iterated
product of order $k$ on $\Sym$ and $\Delta_k$ its iterated coproduct.
These are again the classical Adams operations on $\Sym$, but they are not
algebra morphisms, due to the noncommutativity of $\Sym$: with Sweedler's
notation for the coproduct ($\Delta(F)=F_{(1)}\otimes F_{(2)}$),
\begin{equation}
\begin{split}
\psi^2(FG)=&
\mu\circ\Delta (FG)=\sum_{(F),(G)}F_{(1)}G_{(1)}F_{(2)}G_{(2)}\\
&\not =
\sum_{(F),(G)}F_{(1)}F_{(2)}G_{(1)}G_{(2)} =\psi^2(F)\psi^2(G).
\end{split}
\end{equation}
Similarly to what happens on $QSym$, the $\psi^k$ are given by left internal
product with the reproducing kernel $\sigma_1(kA)$. However, one must pay
attention to the fact that there is another family of such operations,
corresponding to the right internal product with $\sigma_1(kA)$:
\begin{equation}
\psi^k(F(A))= \sigma_1(kA) * F(A) \not = F(A) * \sigma_1(kA) = F(kA).
\end{equation}
The right internal product with $\sigma_1(kA)$ is an algebra morphism (this
follows from the splitting formula~(\ref{splitF})); in terms of alphabets,
this operation corresponds to the transformation $F(A)\mapsto F(kA)$ -we refer
to \cite{ncsf2} for a detailed study of transformations of alphabets in the
framework of noncommutative symmetric functions.
Since they are associated with the same kernels $\sigma_1(kA)$, the spectral
projectors of both families are encoded by the same noncommutative symmetric
functions, the only difference being that one has to take internal product on
different sides.

Thus, the adjoint of our quasi-Eulerian idempotent $e_1$ acting on $QSym$ is
$F\mapsto E_1 * F$ where  $E_1 = \Phi(1) = \log\sigma_1$ is the usual Eulerian
idempotent.

Recall from \cite{GKLLRT} that its action on a product of \emph{primitive}
elements $F_1\cdots F_r$ is given by
\begin{equation}
E_1*(F_1\cdots F_r) = (F_1F_2\cdots)\cdot E_1\,,
\end{equation}
where, on the right hand-side, $E_1$ is the Eulerian idempotent viewed as an element of the group algebra of the symmetric group of order $r$ acting
by \emph{permutation of the indices, e.g.},
$E_1*(F_1F_2)=(1/2)(F_1F_2-F_2F_1)$.

Let us now choose a basis $Q_L$ of the primitive Lie algebra $\mathcal{L}$
of $\Sym$. For the sake of definiteness, we may choose the
Lyndon basis on the sequence of generators $\Phi_n$ (see \cite{GKLLRT}) and we
may assume that $L$ runs over Lyndon compositions.
We can then extend it to a Poincar\'e-Birkhoff-Witt basis of its universal
enveloping algebra $U(\mathcal{L})=\Sym$, so that the Eulerian idempotent will
act by $E_1*Q_I = 0$ if $I$ is not Lyndon, and  $=Q_I$ otherwise \cite{sol1}.
Let now $P_I$ be the dual basis of $Q_I$ in $QSym$. Then, $e_1$ acts by
$e_1(P_I)=0$ if $I$ is not Lyndon, and $=P_I$ otherwise.
Then, $QSym$ is free as a polynomial algebra over the  $P_L$ by Radford's
theorem \cite{Re}. Hence, $e_1$ maps any basis of $QSym$ to a generating set.
Moreover, with our particular choice of the basis, $S^I$ is triangular on the
$Q_I$ so that $P_I$ is triangular on the $M_I$, thus $e_1(M_L)$ for $L$ Lyndon
form a free generating set.

We shall see in the next section that it is true in general that $e_1$
projects the quasi-shuffle algebra onto a generating subspace, although, as we
already mentioned, one can not use any more in the general situation duality
together with the properties of the Eulerian idempotents acting on envelopping
algebras. However, the case of $QSym$ is essentially generic for a wide class
of quasi-shuffle algebras.

A possible line of argumentation (to be developed in a subsequent paper) would
be as follows.
It has been observed in~\cite{CHNT} that noncommutative symmetric functions
provided a good framework for understanding Ecalle's formalism of moulds in a
special case. To deal with the general case, one can introduce the following
straightforward generalization of noncommutative symmetric functions, which
will also provide us with a better understanding of quasi-shuffle algebras in
general.

Let $\Omega$ be an additive monoid, such that any element $\omega$ has only a
finite number of decompositions\footnote{This restriction is not strictly
necessary, but relaxing it would require a generalization of the notion of
Hopf algebra.} 
$\omega =\alpha+\beta$. Let $\Sym^\Omega$ be the free associative algebra over
indeterminates $S_\omega$ ($\omega\in\Omega$, $S_0=1$), graded by
$\deg(S_\omega=\omega$), endowed with the coproduct
\begin{equation}
\Delta S_\omega=\sum_{\alpha+\beta=\omega}S_\alpha\otimes S_\beta
\end{equation}
and define $QSym^\Omega$ as its graded dual.
Let $M_{\bm{\omega}}$ be the dual basis of  $S^{\bm{\omega}}$
Its multiplication rule is obviously given by the quasi-shuffle
(over the algebra $\K[\Omega]$). Let
$\Phi=\log\sum_\omega S_\omega=\sum_\omega \Phi_\omega$.
From this, we can build a basis $\Phi^{\bm{\omega}}$,
a basis of products of primitive elements, which multiplies
by the ordinary shuffle product over the alphabet $\Omega$.

This proves in particular that the quasi-shuffle algebra is isomophic to the
shuffle algebra (here it is one and the same algebra, seen in two different
bases).

Now, the previous argumentation could be copied verbatim here, replacing
compositions by words over $\Omega-\{0\}$.
One may expect that $\Sym^\Omega$ admits an internal product, allowing
to reproduce Ecalle's mould composition. This is the case for example
when $\Omega=\N^r$ \cite{NTcolor}, the algebra $\Sym^\Omega$ being
in this case the natural noncommutative version of McMahon's
multisymmetric functions.

\section{The generalized Eulerian idempotent as a canonical projection}

The purpose of this last section is to show that the results in the previous
section can be, to a large extent, generalized to arbitrary quasi-shuffle
algebras.

It is well-known that quasi-shuffle algebras $QS(A)$ are free commutative
algebras over $A$ (viewed as a vector space). This follows from the
observation that the highest degree component of the quasi-shuffle product map
$QS(A)_n\otimes QS(A)_m\longrightarrow QS(A)_{n+m}\subset QS(A)$ is simply the
usual shuffle product, so that the freeness of $QS(A)$ follows from the
freeness of the tensor algebra over $A$ equipped with the shuffle product by a
standard triangularity argument.

For classical graded commutative connected Hopf algebras $H$ (such as the
tensor algebra over $A$ equipped with the shuffle product), the Leray theorem
asserts that $H$ is a free commutative algebra and one can compute by purely
combinatorial means a family of generators of $H$ as a polynomial algebra
\cite{Pat92,Pat93}. The same result actually holds for quasi-shuffle algebras.
The section is devoted to its proof.

Recall that the characteristic subalgebra $\mathbf{Car}$ of $\WQSym$ is the
free associative subalgebra generated by $I$ (the identity map, when elements
of $\WQSym$ are viewed as operations on quasi-shuffle algebras). This element
$I$ is grouplike in $\WQSym$ and the generalized Eulerian idempotent
$e_1=\log(I)$ is therefore a primitive element. The following proposition
shows that the Hopf algebra structure inherited by $\mathbf{Car}$ from
$\WQSym$ is actually compatible with its action on quasi-shuffle algebras. We
use the Sweedler notation: $\Delta(\sigma)=\sigma^{(1)}\otimes\sigma^{(2)}$.

\begin{proposition}
Let $\sigma\in \mathbf{Car}$, then, we have, for an arbitrary commutative
algebra $A$:
\begin{equation}
((a_1\otimes \dots\otimes a_n)
 \uplus (b_1\otimes \dots\otimes b_m))\cdot \sigma
 = (a_1\otimes\dots\otimes a_n)
   \cdot \sigma^{(1)} \uplus(b_1\otimes \dots\otimes b_m))\cdot\sigma^{(2)}.
\end{equation}
\end{proposition}

The Proposition is obviously true when $\sigma = I$. From
\cite[Lemma~3.1]{egp07}, it is also true in the convolution algebra generated
by $I$, from which the Proposition follows.

Notice that this property is not true for $\WQSym$ -the coproduct in $\WQSym$
is not compatible with the action on quasi-shuffle algebras. This property is
actually already true for $\FQSym$: the coproduct of $\FQSym$ is not
compatible, in general, with the Hopf algebra structure of shuffle or tensors
algebras (this corresponds, in terms of quasi-shuffle algebras, to the
particular case of commutative algebras with a null product).

\begin{theorem}
The generalized Eulerian idempotent $e_1$ is a projection onto a vector space
generating $QS(A)$ as a free commutative algebra (that is, equivalently, the
image of $e_1$ is naturally isomorphic to the indecomposables of $QS(A)$: the
quotient of the augmentation ideal by its square, which is spanned by non
trivial products).
\end{theorem}

We already know that $e_1$ is a spectral idempotent (it maps $QS(A)$ to the
eigenspace of the Adams operations associated with their lowest nontrivial
eigenvalue). From the previous proposition, we have, since $e_1$ is primitive,
for $a,b$ elements of $QS(A)_n,QS(A)_m,\ n,m\not=0$:
\begin{equation}
(a\uplus b)\cdot e_1
=(a\cdot e_1)\uplus (b\cdot \epsilon)+(a\cdot \epsilon)\uplus (b\cdot e_1),
\end{equation}
where $\epsilon$ is the augmentation of $QS(A)$, so that
$b\cdot \epsilon =a\cdot \epsilon=0$, and finally:
\begin{equation}
(a\uplus b)\cdot e_1=0
\end{equation}
from which it follows that the kernel of $e_1$ contains the square of the
augmentation ideal $QS(A)^+$. In particular, $e_1$ induces a well-defined map
on the indecomposables ${\rm Ind} :=QS(A)^+/(QS(A)^+)^2$.

Let us show finally that this map is the identity, from which the Theorem
will follow.
Let us compute first $\Psi^2(a)$ for an arbitrary $a\in QS(A)_n^+$. We have:
\begin{equation}
\Psi^2(a)=2e_1(a)+...+2^ne_n(a),
\end{equation}
where the $e_i$ are higher convolution powers of $e_1=\log(I)$ and map
therefore into $(QS(A)^+)^2$, so that on the indecomposables, $\Psi^2=2e_1$.
On the other hand (using the Sweedler notation),
\begin{equation}
\Psi^2(a)=I^{\ast 2}(a)=2a+a^{(1)}a^{(2)},
\end{equation}
so that, on Ind, $\Psi^2=2 \cdot I$ and the Theorem follows.

\end{document}